\documentstyle[a4,12pt]{article}
\normalbaselines

\newcommand{\spar} {\vskip 0.25cm}

\newcommand{\slp} {straight--line program}
\newcommand{\slps} {straight--line programs}

\newtheorem{proposition}{Proposition}

\newtheorem{theorem}[proposition]{Theorem}

\def\ifm#1#2{\relax\ifmmode#1\else#2\fi}

\newcommand{\klk}    {\ifm {,\ldots,} {$,\ldots,$}}

\def \N{{\rm I\kern -2.1pt N\hskip 1pt}} 
\def \Z{{\ifm{{\rm Z\!\!Z}}{${\rm Z\!\!Z}$}}}
\def \R{{\rm I\kern -2.2pt R\hskip 1pt}} 
\def \K {{\rm I\kern -2.1pt K\hskip 1pt}} 
\def \F {{\rm I\kern -2.1pt F\hskip 1pt}} 
\def \P{{\rm I\kern -2.2pt P\hskip 1pt}}

\newcommand{\Barq}  {\vrule height0.65em width0.05em depth0em \,}
\newcommand{\Q}        {\ifm {\mbox{\rm Q}\hspace{-0.54em}\Barq\>\;}
                             {$\mbox{\rm Q}\hspace{-0.54em}\Barq\>\;$}}
\newcommand{\C}        {\ifm {\mbox{\rm C}\hspace{-0.45em}\Barq\>}
                             {$\mbox{\rm C}\hspace{-0.45em}\Barq\>$}}

\newcommand{\om}[2]   {{#1}_1 \klk {#1}_{#2}} 

\newcommand{\xo}[1]  {\ifm {\om X {#1}} {$\om X {#1}$}}
\newcommand{\xon}    {\ifm {\om X n} {$\om X n$}}

\newcommand{\To}[1] {\ifm {\om T {#1}} {$\om T {#1}$}} 
 
\newcommand{\Uo}[1] {\ifm {\om U {#1}} {$\om U {#1}$}} 
 
\newcommand{\Go}[1] {\ifm {\om G {#1}} {$\om G {#1}$}} 
\newcommand{\Gon} {\ifm {\om G n} {$\om G n$}} 
\def \A {{ {\it I} \! \kern -2.2pt {\rm A}\hskip 1pt}}

\begin{document}
\title{The intrinsic complexity of parametric elimination methods}
\author {J. Heintz, G. Matera, L. M. Pardo, R. Wachenchauzer}
\maketitle
\begin{abstract}
This paper is devoted to the complexity analysis of a particular property,
called {\em algebraic robustness} owned by all known symbolic methods of 
parametric polynomial equation solving (geometric elimination).
It is shown that {\em any} parametric elimination procedure which owns this
property must neccessarily have an exponential sequential time complexity 
even if highly performant data structures 
(as e.g. the \slp\ encoding of polynomials) are used.
The paper finishes with the motivated introduction of a new non-uniform
complexity measure for zero-dimensional polynomial equation systems, called
{\em elimination complexity}.
\end{abstract}
\smallskip


\bigskip

\centerline{Research was partially supported by the following
Argentinian and Spanish grants~:}

\centerline{UBA-CYT.EX.001, PIP CONICET 4571, DGICYT PB93--0472--C02--02.}
 

\begin{section}{Introduction}\label{sec-intro}
Modern algebraic geometry startd about 200 years ago as {\em algorithmic}
algebraic geometry, and, more precisely, as algorithmic elimination theory.
The motivation for the creation of such a field was the search for methods
which allow to find the {\em real} solutions of a polynomial equation system.
Nevertheless, the very origin of algebraic geometry was given by the
observation that real root finding is a rather unfeasable task without a
previous study of the behaviour of the {\em complex}
solutions of polynomial equation systems. This observation was first made by
Euler and B\'ezout and then extended to a general theory by a long list of
geometers of the last century.
This list includes names as Jacobi, Sylvester, Kronecker, M. Noether,
Hilbert (the creator of modern commutative algebra), Castelnuovo, Bertini,
Enriques.
\spar

Unfortunately modern algebraic geometry deviated from its origins
(mainly after World War II) and todays written contributions to this field
are authentic masterpieces in the art of boring the reader.
Nevertheless in the last twenty years a new community of old algebraic
geometers doing symbolic computation splitted out of the mainstream. The
intention of this community to bring back algebraic geometry to its origin
(and to introduce also new aspects like efficient polynomial equation
solving for industrial applications) must be praised highly.
On the other hand the (mainly rewriting based) computational approach used
by this community is far too simple minded for the difficult task of
efficient, i.e. real world polynomial equation solving. An important
drawback of this approach consists in the almost total absence of todays
skill in algorithmics and data structure manipulation as well as the
unawareness of modern programming techniques coming from software
engineering. There is no place here to describe in detail the advances and
weaknesses of symbolic computation techniques applied to elimination theory
(computer algebra). For an overview about rewriting based methods (Groebner
basis techniques) we refer to the books \cite{BeWe93}, 
\cite{Mishra93}, \cite{CoLiOS92}
(these books include also motivations and historical considerations).
The state of the art in sparse techniques can be found in \cite{Emiris96}.
Finally the seminumerical approach to elimination theory is described in 
the book \cite{BlCuShSm97} and the surveys
\cite{HeMo93}, \cite{Pardo95}, \cite{Heintz97b} and in the
research papers \cite{GiHeMoMoPa95}, \cite{GiHaHeMoMoPa96}
and \cite{BaGiHeMaMb97}.
\spar


Informally we call
a parametric elimination procedure {\em algebraically robust} if it produces
for (geometrically or scheme-theoretically) flat families of problem
instances ``continuous'' or ``stable'' solutions.
\spar

Of course this notion of algebraic robustness depends on the (geometric or
scheme-theoretical) context, i.e. it is not the same for schemes or
varieties. Below we are going to explain our idea of algebraic robustness in the typical
situation of flat families of algebraic varieties given by reduced complete
intersections.
\spar

The main outcome of this paper is the observation that all known parametric
elimination procedures are algebraically robust and that all algebraically
robust elimination procedures must neccessarily have an exponential time
complexity even if the highly performant encoding of polynomials by \slps\
is used (Theorem 1).
Therefore a revolutionary change of mathematical theory and algorithmics would
be neccessary in order to design a (possibly non-existent) highly performant
general purpose elimination procedure.
The rest of the paper is devoted to the motivated introduction of a new
uniform complexity measure
for zero-dimensional polynomial equation systems, called
{\em elimination complexity}.


\begin{subsection}{Flat families of elimination problems}\label{subsec11}

Let $k$ be an infinite and perfect field
with algebraic closure $\bar{k}$ and let
$\To m,$
$\Uo r, \xon, Y$ be indeterminates over $k$. 
Let $\Gon, F$ be
polynomials belonging to the $k$-algebra $k[\To m, \Uo r, \xon].$
Suppose that the polynomials $\Go n$ form a regular sequence in $k[\To m,
\Uo r, \xon]$ defining thus an equidimensional subvariety $V:=\{ G_1=0\klk
G_n=0\}$ of 
the $(m+r+n)$-dimensional affine space
$\A^{m+r+n}$ over the field $\bar{k}$. 
The algebraic variety $V$ has dimension $m+r$.
Let $\delta$ be the
(geometric) degree of $V$. 
Suppose furthermore that the morphism of affine varieties
$\pi:V\longrightarrow \A^{m+r}$, induced by the canonical projection of
$ \A^{m+r+n}$ onto $ \A^{m+r}$, is finite and generically unramified (this
implies that $\pi$ is flat).
Let $\tilde{\pi}:V\longrightarrow \A^{m+r+1}$ be the morphism defined by
$\tilde{\pi}(z) := (\pi(z),F(z))$ for any point $z$ of the variety $V$.
The image of $\tilde{\pi}$ is a hypersurface of $\A^{m+r+1}$ whose minimal
equation is a polynomial
of $ k[\To m, \Uo r, Y]$ which we denote by $P$. Observe that $P$ is monic
in $Y$ and that $\deg P\leq \delta$ holds. Furthermore $\deg _Y P$ is the
cardinality of the image of the restriction of $F$ to the set $\{w\}\times
\pi^{-1}(w)$, where $w$ is a typical point of $\A^{m+r}$.
The polynomial $P(\To m,\Uo r,F)$ vanishes on the variety $V$. \spar


Let us consider an arbitrary point $t=(t_1,\ldots ,t_m)$ of $\A ^m$. For
arbitrary polynomials $B\in k[\To m,\Uo r,\xon]$ and $C\in k[\To m,\Uo
r,Y]$ we denote
by $B^{(t)}$ and $C^{(t)}$ the polynomials $B(t_1\klk t_m,\Uo r,\xon)$ and
$C(t_1\klk t_m,\Uo r,Y)$ which belong to $k(t_1\klk t_m)[\Uo r,\xon]$ and
$k(t_1\klk t_m)[\Uo r,Y]$ respectively.
Similarly we denote for an arbitrary polynomial $A\in k[T_1,\ldots, T_m]$
by $A^{(t)}$ the value $A(t_1,\ldots ,t_m)$ which belongs to the field
$k(t_1,\ldots ,t_m)$. The polynomials $G_1^{(t)}\klk G_n^{(t)}$ form a
regular sequence in
$k(t_1\klk t_m)[\Uo r,\xon]$ and define an equidimensional subvariety $V^{(t)}:=
\{G_1^{(t)}=0\klk G_n^{(t)}=0\}$ of $\A^{r+n}$ whose degree is bounded by
$\delta$.
Let $\pi^{(t)}:\ V^{(t)} \longrightarrow \A^r$ and $\tilde{\pi}^{(t)}:\
V^{(t)} \longrightarrow \A^{r+1}$ be the morphisms induced by $\pi$ and
$\tilde{\pi}$ on the variety $V^{(t)}$.
Then the morphism $\pi^{(t)}$ is finite and flat but not necessarily
generically unramified. 
Furthermore the image of $\tilde{\pi}^{(t)}$ is a hypersurface of
$\A^{r+1}$ on which the polynomial $P^{(t)}$ vanishes (however
$P^{(t)}$ is not necessarily the minimal equation of this hypersurface). We
call the equation system
$G_1 =0\klk G_n =0$ and the polynomial $F$ {\sl a flat family of
$r$-dimensional elimination problems depending on the parameters $\To m$}.
An element $t\in\A^m$ is considered as a {\sl parameter point} which
determines a {\sl particular problem instance}. The equation system $G_1
=0\klk G_n =0$ together with the polynomial $F$ is called the {\sl general
instance} of the given elimination problem and the polynomial $P$ is called
its {\sl general solution}.\\


The problem instance determined by the parameter point $t\in \A^m$ is given
by the equations $G_1^{(t)}=0\klk G_n^{(t)}=0$ and the polynomial
$F^{(t)}$. The polynomial $P^{(t)}$ is called the {\em solution} of this
particular problem instance. We call two parameter points $t,t' \in \A^m$
{\sl equivalent} (in symbols $t \sim t'$) if $G_1^{(t)}=G_1^{(t')} \klk
G_n^{(t)}=G_n^{(t')}$ and $F^{(t)} = F^{(t')}$ holds.
Observe that $t \sim t'$ implies also $P^{(t)} = P^{(t')}$. We call
polynomials $A\in k[\To m]$,
$B\in k[\To m,\Uo r,\xon]$ and $C\in k[\To m,\Uo r,Y]$ invariant (with
respect to $\sim$) if for any two parameter points $t,t'$ of $\A^m$ with
$t\sim t'$ the identities $A^{(t)} = A^{(t')}$, $B^{(t)} = B^{(t')}$ and
$C^{(t)}=C^{(t')}$ hold. \spar


A straight-line program in
$k[T_1,\ldots ,T_m,U_1,\ldots ,U_r,Y]$ {\sl with parameters in
$k[T_1,\ldots ,T_m]$} is an arithmetic circuit in $k[T_1,\ldots
,T_m,U_1,\ldots ,U_r,Y]$, say $\gamma$, modelized in the following way:
$\gamma$ is given by a directed acyclic graph whose internal nodes are
labeled as usually by arithmetic operations. The input nodes of $\gamma$
are labeled by the variables $U_1,\ldots ,U_r$ and the nodes of $\gamma$ of
indegree $0$ which are not input nodes (i.e. the parameter nodes of
$\gamma$) are labeled by arbitrary elements of
$k[T_1,\ldots ,T_m]$ called {\sl parameters of $\gamma$}. We call such a
straight-line program
$\gamma$ {\sl invariant } (with respect to the equivalence relation $\sim
$) if all its parameters are invariant polynomials of
$k[T_1,\ldots ,T_m]$. Let us observe that in an absolutely analogous way
one may extend the notion of a straight-line program with parameters in
$k[T_1,\ldots ,T_m]$
and the notion of an invariant straight-line program to arithmetic circuits
defined in
$k[T_1,\ldots ,T_m,U_1,\ldots ,U_r,X_1,\ldots ,X_n]$ or in $k[T_1,\ldots
,T_m,U_1,\ldots ,U_r]$ (however we shall make almost exclusive use of the
corresponding notion for circuits in\\ 
$k[T_1,\ldots ,T_m,U_1,\ldots,U_r,Y]$). \spar

We are now ready to characterize in the given situation what we mean by an
{\sl algebraically robust elimination procedure}. Suppose that the
polynomials $\Go n$ and $F$ are given by a \slp $\beta$ in $k[\To m,\Uo
r,\xon]$. An algebraically robust elimination procedure produces from the
circuit $\beta$ as input an {\sl invariant} \slp $\Gamma$ in $k[\To m,\Uo
r, Y]$ as output such that $\Gamma$ represents the polynomial $P$. Observe
that in our definition
of algebraic robustness we did not require that $\beta$ is an invariant
\slp because
this would be too restrictive
for the modeling of real
situations in computational
elimination theory. However in
all known practical situations $\beta$ always satisfies the following (weak)
invariance condition: any specialization of the input variables $\xo n $
in $\beta$ into arbitrary values of $k$ transforms the circuit $\beta$ in a
invariant \slp\ in $k[\To m, \Uo r ]$.
\spar

The invariance property required for the output circuit $\Gamma$ means the
following: let $t=(t_1\klk t_m)$ be a parameter point of $\A^m$ and let
$\Gamma^{(t)}$ be the \slp in $k(t_1\klk t_m)[\Uo r, Y]$ obtained from the
circuit $\Gamma$ evaluating in $t$ the elements of $k[\To m]$ which occur
as parameters of $\Gamma$. Then the \slp $\Gamma^{(t)}$ depends only on the
particular problem instance determined by the parameter point $t$ but not
on $t$ itself. Said otherwise, an algebraically robust elimination
procedure produces the solution of a particular problem instance in a way
which is independent of the possibly different parameter points which
may determine the given problem instance.
\spar


Observe that by definition an algebraically robust elimination procedure
produces always the general solution of the elimination problem under
consideration. We are going to show - informally - that none of the known
(exponential time) parametric elimination procedures can be transformed into
a polynomial time algorithm. For this purpose it is important to remark that
the known parametric elimination procedures (which are without exception
based on linear algebra as well as on comprehensive Groebner basis
techniques) are all algebraically robust on flat families of elimination
problems if the general solution of the given problem is required as output.


The invariance property of these procedures is easily verified in the
situation of the flat family of $r$-dimensional elimination problems
introduced before. One has only to observe that all known elimination
procedures accept the input polynomials $\Gon, F$ in their dense or sparse
coefficient representation or as evaluation black box with respect to the
variables $\Uo r,\xon$.

\end{subsection} 


\begin{subsection}{A particular flat family of 1-dimensional elimination
problems}

Let $S, T, U, \xo {2n} , Y$ be indeterminates over $\Q$. We consider the
following flat family of one-dimensional elimination problems depending on
the parameters $S$ and $T$. Let

$$G_1:=X_1^2-X_1\klk G_n:=X_n^2-X_n,$$

$$G_{n+1} := X_{n+1}-\sum_{1\leq k\leq n} 2^{k-1}X_k,
G_{n+2}:=X_{n+2}-X_{n+1}^2, G_{2n}:=X_{2n}-X_{2n-1}^2$$

and

$$F:=(1+S\prod_{1\leq i\leq n} (T^{2^{i-1}}+X_{n+i})) \prod_{1\leq j\leq n}
((U^{2^{j-1}}-1)X_j+1).$$

We interprete $G_1\klk G_{2n}$ and $F$ as elements of the polynomial ring\\
$\Q[S,T,U,\xo {2n}]$. Thus we have $m:=2$ and $r:=1$ in this situation.


It is clear from this representation that the polynomials $\Go {2n}$ and
$F$ can be evaluated by a \slp $\beta$ in $\Q[S,T,U, \xo{2n}]$ of nonscalar
length $O(n)$ (which satisfies the invariance condition for input circuits
mentioned before). The degree in $\xo {2n}$ of the polynomials $\Go {2n}$
and $F$ is bounded by $2n$. The variety $V:= \{ G_1=0\klk G_n = 0\}$ is the
union of $2^n$ affine linear subspaces of $\A^{3+2n}$ of the form
$\A^3\times \{ \xi\}$, where $\xi$ is a solution of the equation system
$G_1 = 0\klk G_{2n}=0$ in $\Z^{2n}$. The morphism $\pi:\ V \rightarrow
\A^3$ is obtained by gluing together the canonical projections onto $\A^3$
of these affine linear spaces. Obviously the morphism$\pi$ is finite an
generically unramified. In particular $\pi$ has constant fibers. Let
$(l_1\klk l_n)$ be a point of $\{0,1\}^n$ and let $l:=\sum_{1\leq j\leq n}
l_j 2^{j-2}$ be the integer $0\leq l< 2^n$ with bit representation $l_n
l_{n-1}\ldots l_1$.
Put $l_{n+1}:=l,l_{n+2}:=l^2\klk l_{2n}:=l^{2^{n-1}}$. One verifies immediately
that with the convention $0^0: = 1$ the identity

$$F(S,T,U,l_1\klk l_{2n}) = U^l(1+S\sum_{ {0\leq p,q<2^n}\atop {p+q=2^n-1}}
T^p l^q)$$

holds. Therefore for any point $(s,t,u,l_1\klk l_{2n})\in V$ with
$l:=\sum_{1\leq j\leq n} l_j 2^{j-1}$ we have

$$F(s,t,u,l_1\klk l_{2n}) = u^l(1+s\sum_{{0\leq p,q<2^n}\atop {p+q=2^n-1}}
t^p l^q).$$

From this we deduce easily that the elimination polynomial
$P\in \Q[S,T,U,Y]$ we are looking for is


$$P:= \prod_{0\leq l<{2^n}} \left( Y - U^l(1+S\sum_{{0\leq p,q<2^n}\atop
{p+q=2^n-1}} T^p l^q)\right).$$

This polynomial has the form

$$P= Y^{2^n} - \left( \sum_{0\leq l<{2^n}} U^l(1+S\sum_{{0\leq
p,q<2^n}\atop{ p+q=2^n-1}} T^p l^q)\right)Y^{2^n-1} + \hbox{ lower degree
terms in }Y.$$

Suppose now that there is given an algebraically robust elimination
procedure. This procedure produces from the input circuit $\beta$ an
invariant \slp $\Gamma$ in $\Q[S,T,U,Y]$, which evaluates the polynomial
$P$. Recall that the invariance of $\Gamma$ means that the parameters of
$\Gamma$ are invariant polynomials of $\Q[S,T]$, say $A_1\klk A_N$. \spar


Let ${\cal L}(\Gamma)$ be the total and $L(\Gamma)$ the non-scalar length
of the \slp $\Gamma$. We have $L(\Gamma) \leq {\cal L}(\Gamma)$ and $N\leq
(L(\Gamma)+3)^2$. Let $Z_1\klk Z_N$ be new indeterminates. From the graph
structure of the circuit $\Gamma$ we deduce that there exist for $0\leq
l<2^n$ polynomials $Q_l\in\Z[Z_1 \klk Z_N]$ such that $Q_l(A_1\klk A_N)$ is
the coefficient of the monomial $U^l Y^{2^n-1}$ of $P$. This means that we
have for $0\leq l<2^n$ the identity

\begin{equation}\label{eq1}
Q_l(A_1\klk A_N) = 1+S\sum_{{0\leq p,q<2^n}\atop{p+q=2^n-1}} T^p l^q.
\end{equation}


Observe now that for any two values $t,t'\in\A^1$ the points $(0,t)$ and
$(0,t')$ of $\A^2$ are equivalent (in symbols: $(0,t)\sim (0,t')$). From
the invariance of $A_1\klk A_N$ we deduce therefore that $A_j(0,t) =
A_j(0,t')$ holds for $1\leq j\leq N$. This means that $\alpha_1 :=
A_1(0,T)\klk \alpha_N := A_N(0,T)$ are constant values of $\Q$. From 
identity (\ref{eq1}) we deduce that $Q_l(\alpha_1\klk \alpha_N) = 1$
holds for any $0\leq l <2^n$.

\spar

Let us consider the morphisms of affine spaces $\mu:\ \A^2 \longrightarrow
\A^N$ and $\psi:\ \A^N\longrightarrow \A^{2^n}$ given by $\mu:=(A_1\klk
A_N)$ and $\psi:=(Q_l)_{0\leq l<{2^n}}$. \spar Observe that

$$\psi\circ \mu = ( Q_l(A_1\klk A_N))_{0\leq l<{2^n}} = (1+S\sum_{{0\leq
p,q<2^n}\atop{ p+q={2^n}-1}} T^p l^q)_{0\leq l<{2^n}} $$

holds. Furthermore let us consider the points $\alpha := (\alpha_1\klk
\alpha_N)\in \A^N$ and $\omega :=(1\klk 1)\in \A^{2^n}$. From 
our previous considerations we deduce the identities

\begin{eqnarray*}
(\psi\circ \mu) (0,T) = & (Q_l(A_1(0,T)\klk A_N(0,T)))_{0\leq l<{2^n}} = & \\
= & ( Q_l(\alpha_1\klk\alpha_N))_{0\leq l <{2^n}} = (Q_l(\alpha))_{0\leq
l<{2^n}} & = (1\klk 1) = \omega.
\end{eqnarray*}

In particular we have $\psi(\alpha) = \omega$. We analyze now 
the local behaviour of 
the morphism
$\psi$ in the point $\alpha \in\A^N$. Let ${E_{\alpha}}$ and
${E_{\omega}}$ be the tangent spaces of
the points $\alpha$ and $\omega$ belonging to $\A^N$ and $\A^{2^n}$
respectively. Let
us denote the differential of the map $\psi$ in the point $\alpha$ by $(D
{\psi})_{\alpha}:\ {E_{\alpha}}\longrightarrow {E_{\omega}}$. Taking the
canonical projections of $\A^N$ and $\A^{2^n}$ as local coordinates in the
points $\alpha$ and $\omega$ respectively, we identify ${E_{\alpha}}$ with
$\A^N$ and ${E_{\omega}}$ with $\A^{2^n}$.

For any value $t\in \Q$ we consider the parametric curves $\gamma_t:\ \A^1
\longrightarrow \A^N$ and $\delta_t:\ \A^1 \longrightarrow \A^{2^n}$
defined by

$$\gamma_t := (A_1(S,t)\klk A_N(S,t))\hbox{ and } \delta_t:=
(1+S\sum_{{0\leq p,q<2^n}\atop{ p+q=2^n-1}} t^p l^q)_{0\leq l<{2^n}}$$

respectively.
Observe that $\psi \circ \gamma_t = \delta_t$ and that $\gamma_t(0) =
\alpha$, $\delta_t(0) = \omega$ holds (independently of the value $t$).
\spar

We consider $\gamma_t$ and $\delta_t$ as one-parameter subgroups of $\A^N$
and $A^{2^n}$ respectively. \spar

Now fix $t\in\Q$ and consider

$$\gamma_t' (0) = ( \frac {\partial A_1} {\partial S } (0,t)\klk \frac
{\partial A_N} {\partial S } (0,t) ) \hbox{ and } \delta_t' (0) =
(\sum_{{0\leq p,q<2^n}\atop { p+q=2^n-1}} t^p l^q)_{0\leq l<{2^n}}.$$


We have $\gamma_t'(0) \in {E_{\alpha}}$ and $\delta'(0) \in {E_{\omega}}$.
Moreover one sees easily that

$$({D}{\psi})_{\alpha}(\gamma_t'(0)) =\delta_t'(0) = ( \sum_{{0\leq
p,q<2^n}\atop{ p+q=2^n-1}} t^p l^q)_{0\leq l<{2^n}}$$

holds. Choosing now $2^n$ different values $t_0\klk t_{2^{n}-1}$ of $\Q$ we
obtain $2^n$ tangent vectors $v_l:=\gamma_{t_l}'(0)$ of $E_{\alpha}$ with
$0\leq l< 2^n$. Let $M$ be the $2^n\times 2^n$ matrix whose row vectors are
$({D}{\psi})_{\alpha}(v_0)\klk
({D}{\psi})_{\alpha}(v_{2^n-1})$.
Observe that $M$ has the form

$$M = (t_{h}^{2^n-p-1})_{0\leq h, p < {2^n}} \cdot (l^q)_{0\leq l,q<{2^n}}.$$


Thus $M$ is the product of two non-singular Vandermonde matrices and
therefore itself non-singular. This means that in ${E_{\omega}}$ the
tangent vectors $( D\psi)_{\alpha}(v_0)\klk (D\psi)_{\alpha}(v_{2^n-1})$
are linearly independent. Hence in ${E_{\alpha}}$ the tangent vectors
$v_0\klk v_{2^n-1}$ are linearly independent too. \spar

Therefore we have $2^n\leq \dim E_{\alpha} = N$ which implies $2^n\leq
N\leq (L(\Gamma)+3)^2$.
>From this we deduce the estimation
$2^{{\frac n 2}}-3 \leq L(\Gamma) \leq {\cal L}(\Gamma)$. \spar

{\em
We have therefore shown that any algebraically robust elimination procedure
applied to our flat family of one-dimensional elimination problems produces
a solution circuit of size at least $2^{{\frac n 2}}-3$ i.e. a circuit of
exponential size in the length $O(n)$ of the input.
}
\spar


The discussion of the previous example
shows that the objective of
a polynomial time procedure for geometric (or algebraic) elimination can
not be reached following a evolutionary way, i.e. constructing improvements
of known elimination methods. \spar

It was fundamental in our argumentation above that our notion of {\em
algebraically robust elimination procedure} excludes branchings in the
output program. This suggests that any polynomial time elimination
algorithm (if there exists one) must have a huge topological complexity.
Thus hypothetical efficiency in geometric elimination seems to imply
complicated casuistics.
\spar


This idea is worth to be discussed further. One may also ask whether
admitting divisions in the output circuit helps to lower its minimal size.
To some limited extent divisions in the output circuit are compatible with
our proof method. However one has to take care of the way how these
divisions may affect the dependence of the coefficients of the output
polynomial on the parameters of the circuit representing it.

\spar

The formulation of a condition which guarantees the generalisation of our
method to output circuits with divisions seems to be cumbersome. In our
example one has to make sure that
even in presence of divisions for any value $t\in\Q$ the one-parameter
subgroup $\gamma_t$ still converges to one and the same point of $\A^N$.

\spar


Finally let us mention that our proof method above contributes absolutely
{\em nothing} to the elucidation of the fundamental thesis of algebraic
complexity theory, which says that geometric elimination is non-polynomial
in the (unrestricted) non-uniform complexity model. Similarly no advance is
obtained by our method with respect to the question whether $P_{\C} \neq
NP_{\C}$\ holds in the BSS complexity model, see \cite[Chapter
7]{BlCuShSm97}. \spar

In fact our contribution consists only in the discovery of a very limiting
uniformity property (algebraic robustness) present in all known elimination
procedures. This uniformity property inhibits the transformation of
these elimination procedures into polynomial time algorithms.
We resume the conclusions from the complexity discussion of our example in
the following form:

\begin{theorem}
For any $n\in\N$ there exists a one-dimensional elimination problem
depending on one parameter and $2n+1$ variables, having input length $O(n)$
such that the following holds:
any algebraically robust elimination procedure which solves this problem
produces an output circuit of size at least $2^{n\over 2}-3$ (i.e. of
exponential size in the input length).
\end{theorem}

\end{subsection}
\end{section}


\begin{section}{On the complexity of geometric elimination procedures in
the unrestricted non-uniform model}

Let us now analyze from a general non-uniform point of view how the
seminumerical elimination procedure
designed in \cite{GiHeMoMoPa95} and \cite{GiHaHeMoMoPa96}
works on a given flat family of zero-dimensional elimination
problems.

Let $\To m,\xon, Y$ be indeterminates over the ground ground field $k$ and
let $\Go n, F$ be polynomials belonging to the $k$-algebra $k[\To m,\xon]$.
Let $d:= \max \{ \deg G_1\klk \deg G_n\}$ and suppose that $\Go n$ and $F$
are given by \slps in $k[\To m\xon]$ of length $L$ and $K$ respectively.
Suppose that the polynomials $ G_1\klk G_n$ form a regular sequence in\\
$k[\To m,\xon]$ defining thus an equidimensional subvariety
$V = \{ G_1= 0\klk G_n=0\}$ of $\A^{m+n}$ of dimension $m$. \spar


Assume that the morphism $\pi:\ V \longrightarrow \A^m$, induced by the
canonical projection of $\A^{m+n}$ onto $\A^m$ is finite and generically
unramified. Let $\delta$ be the degree of the variety $V$ and let $D$ be
the degree of the morphism $\pi$. Furthermore let $\tilde{\pi}:\
V\longrightarrow \A^{m+1}$ be the morphism of affine varieties defined by
$\tilde{\pi}(z) := (\pi(z), F(z))$ for any point $z$ of $V$. Let $P\in
k[\To m,Y]$ be the
minimal polynomial of the image of $\tilde{\pi}$. The polynomial $P$ is
monic in $Y$ and one sees immediatly that $\deg P\leq \delta \deg F$ and
$\deg_Y P \leq D$ holds. Let us write $\delta_{*}:=\deg_{\To m} P$ and
$D_{*}:=\deg_Y P$. \spar


Let us consider as Algorithm 1 and Algorithm 2 two non-uniform variants of
the basic elimination method designed in \cite{GiHeMoMoPa95} and
\cite{GiHaHeMoMoPa96}.

\begin{itemize}

\item Algorithm 1 is represented by an arithmetic network of size
$K\delta^{O(1)}+L( n d \Delta)^{O(1)}$ where $\Delta$ is the degree of the
equation system $G_1=0\klk G_n=0$ (observe that $\Delta\leq \deg G_1\cdots
\deg G_n$ always holds). The output is a \slp $\Gamma_1$ in $k[\To m,Y]$ of
length $(K+L)(n \delta)^{O(1)}$ which represents the polynomial $P$.

\item Algorithm 2 starts from the geometric description of a unramified
parameter (and lifting) point $t = (t_1\klk t_m)$ of $k^m$ which has the
additional property that the image of $F$ restricted to the set
$\{t\}\times \pi^{-1}(t)$ has cardinality $D_{*}=\deg_Y P$. The algorithm
produces then a
\slp $\Gamma_2$ in $k[\To m,Y]$ of length $O(KD_{*}^{O(1)} \log \delta_{*})
+ \delta_{*}^{O(1)} = K(\delta \deg F)^{O(1)}$ which represents the
polynomial $P$.

\end{itemize}


We observe that $K\delta^{O(1)}$ is a characteristic quantity which appears
in the length of both \slps $\Gamma_1$ and $\Gamma_2$. We are going now to
analyze the question whether
a complexity of type $K\delta^{O(1)}$ is intrinsic for the elimination
problem under consideration.

\spar

In the next subsection we are going to exhibit an example of a particular
zero-dimensional elimination problem for which the quantity $K\delta$
appears as a lower bound for the nonscalar size of the output polynomial in
the {\em unrestricted non-uniform} complexity model.



\begin{subsection}{A particular flat family of zero-dimensional elimination
problems}

Let $S,\To {\delta}, X,Y$ be indeterminates over $\Q$ and let\\
$G:=\prod_{1\leq l\leq \delta}(X - T_l)$ and $F:=SX^{2^K}$. Let
$V:=\{G=0\}$ be the hypersurface of $\A^{\delta+2}$ defined by the
polynomial $G$ and let $\pi:\ V\longrightarrow \A^{\delta+1}$ be the
finite, generically unramified morphism induced by the canonical projection
of $\A^{\delta+2}$ onto $\A^{\delta+1}$.

Observe that $\delta$ is the degree of the hypersurface $V$ of
$\A^{\delta+2}$ and of the morphism $\pi$ (in fact $V$ is the union of
$\delta$ distinct hyperplanes of $\A^{\delta+2}$).
\spar


The polynomials $G$ and $F$ have a nonscalar complexity $\delta$ and $K+1$
respectively. They represent a flat family of zero-dimensional elimination
problems with $m:=\delta+1, n:=1, \deg V = \delta$ and $\deg \pi = \delta$
(see Subsection \ref{subsec11}).
The general solution of this elimination problem is represented by the
polynomial

$$P:=\prod_{1\leq l\leq \delta} (Y-ST_l^{2^K}) =
Y^{\delta}-Y^{\delta-1}S\sum_{1\leq l \leq \delta } T_l^{2^K}+\hbox{ higher
degree terms in }S$$

which belongs to $\Q[S,\To {\delta},Y]$. Let $\Gamma$ be s \slp of
nonscalar length $L(\Gamma)$ in $k[S,\To {\delta},Y]$ which computes the
polynomial $P$. We transform the circuit $\Gamma$ into a \slp $\Gamma^*$ in
$\Q[\To{\delta}]$ of nonscalar length $L(\Gamma^*)\leq 3 L(\Gamma)$ which
computes the polynomial $R:=\sum_{1\leq l\leq \delta} T_l^{2^K}$.

This can be done as follows: first we derive the \slp $\Gamma$ with respect
to the variable $S$ and then we specialize $S$ and the variable $Y$ into
the values $0$ and $1$ respectively.
\spar

Analyzing the complexity of the polynomial $R$ by means of Strassen's
degree method as in \cite{BaSt83}, \cite{Strassen73A} we find that
$L(\Gamma^*)\geq (K-1)\delta$ holds. This implies $L(\Gamma)\geq \frac 1 3
(K-1) \delta = \Omega(K\delta)$.

\spar

Unfortunately the meaning of the parameter $\delta$ is ambiguous: $\delta$
is the degree of the variety $V$ and the degree of the morphism $\pi$ as
well as the nonscalar complexity of the polynomial $G$. Nevertheless our
example shows that any
optimal elimination procedure which produces the general solution of a
given flat family of zero-dimensional elimination problems has an inherent
complexity which depends linearly on the nonscalar length of the polynomial
which defines the projection we are considering. The factor of
proportionality of this linear dependence appears as an invariant of the
equational part of our elimination problem. For the moment we are not able
to interpret unambiguously this factor of proportionality. It is always
bounded from above by a polynomial function of the \slp size, of the number
of variables eliminated and of the degree of the input variety and appears
in some cases as bounded from below by a quantity which may be interpreted
alternatively as the degree of the input system or as its nonscalar length.
This leads us to the following notion of {\em elimination complexity} of a
given flat family of zero-dimensional elimination problems. This notion is
the subject of the next
subsection.


\end{subsection} 


\begin{subsection}{The elimination complexity of a zero-dimensional
polynomial equation system}

Let $\To m,\xon,Y$ be indeterminates over $k$ and let \\
$\Go n\in k[\To m,\xon]$ 
be polynomials forming a regular sequence in $k[\To m,\xon]$. Let
$V=\{ G_1 =0\klk G_n = 0\}$ be the equidimensional variety defined by the
polynomials $\Go n$ and let $\pi:\ V \longrightarrow \A^m$ be the morphism
of affine varieties induced by the canonical projections of $\A^{m+n}$ onto
$\A^m$. Suppose that $\pi$ is finite and generically unramified and observe
that $V$ has dimension $m$. \spar


For any polynomial $F\in k[\To m,\xon]$ we consider the flat family of
zero-dimensional elimination problems given by the equations $G_1=0\klk
G_n=0$ and the polynomial $F$.
Let $P_F\in k[\To m, Y]$ be the general solution of this problem. Let $L(F)
$ and $L(P_F)$ be the nonscalar complexity of the polynomials $F$ and $P_F$
respectively. Observe that the set of values

$$N_{\Go n} := \{ \frac {L(P_F)}{L(F)} ; F\in K[\To m,\xon] \}$$

is bounded by a quantity which depends polynomially on the degree of $V$
and the number $n$ of variables to be eliminated (compare this with the
length of the \slp $\Gamma_1$ in Algorithm $1$ in this section).

\spar


We define now the supremum $\sup N_{\Go n}$ as the {\em elimination
complexity} of the equation system $G_1=0\klk G_n=0$. In the example of the
previous subsection the nonscalar
\slp length of the equational part of the input system equals the quantity
$\deg V$ and this quantity represents
a lower bound for the elimination complexity of the given equation system.

\end{subsection} 
\end{section} 


{\small

\begin{thebibliography}{10}

\bibitem{BaGiHeMaMb97}
B.~Bank, M.~Giusti, J.~Heintz, R.~Mandel, and G.~Mbakop.
\newblock Polar Varieties and Efficient Real Equation Solving: The Hypersurface
  Case.
\newblock {\em J. of Complexity\/}, {\bf vol.~13}:pp.~5--27, 1997.

\bibitem{BaSt83}
W.~Baur and V.~Strassen.
\newblock The complexity of partial derivatives.
\newblock {\em Theoret. Comp. Sci.\/}, {\bf vol.~22}:pp.~317--330, 1983.

\bibitem{BlCuShSm97}
L.~Blum, F.~Cucker, M.~Shub, and S.~Smale.
\newblock Complexity and Real Computation.
\newblock preprint, 1997.

\bibitem{CoLiOS92}
D.~Cox, J.~Little, and D.~O'Shea.
\newblock {\em Ideals, Varieties, and Algorithms: an introduction to
  computational algebraic geometry and commutative algebra\/}.
\newblock Undergraduate Texts in Mathematics. Springer Verlag, Berlin, 1992.

\bibitem{Emiris96}
I.~Emiris.
\newblock On the Complexity of Sparse Elimination.
\newblock {\em J.\ Complexity\/}, {\bf vol.~12}:pp.~134--166, 1996.

\bibitem{GiHaHeMoMoPa96}
M.~Giusti, K.~H{\"a}gele, J.~Heintz, J.~E. Morais, J.~L. {Monta\~na}, and L.~M.
  Pardo.
\newblock Lower Bounds for Diophantine Approximation.
\newblock In {\em Proceedings of MEGA'96\/}, vol. 117,118, pp. 277--317.
  Journal of Pure and Applied Algebra, 1997.

\bibitem{GiHeMoMoPa95}
M.~Giusti, J.~Heintz, J.~E. Morais, J.~Morgenstern, and L.~M. Pardo.
\newblock Straight--Line Programs In Geometric Elimination Theory.
\newblock {\em to appear in J. of Pure and App. Algebra\/}, pp. 1--46, 1997.

\bibitem{Heintz97b}
J.~Heintz.
\newblock The virtual world as a real trap for the scientist.
\newblock Invited talk at WAIT'97, Buenos Aires, Argentina, 1997.

\bibitem{HeMo93}
J.~Heintz and J.~Morgenstern.
\newblock On the intrinsic complexity of elimination theory.
\newblock {\em J. of Complexity\/}, {\bf vol.~9}:pp.~471--498, 1993.

\bibitem{Mishra93}
B.~Mishra.
\newblock {\em Algorithmic Algebra\/}.
\newblock Springer Verlag, New York, 1993.
\newblock ISBN 0-387-94090-1.

\bibitem{Pardo95}
L.~M. Pardo.
\newblock How lower and upper complexity bounds meet in elimination theory.
\newblock {\em AAECC-11\/}, 1995.

\bibitem{Strassen73A}
V.~Strassen.
\newblock Die Berechnungskomplexit{\"a}t von elementarsymmetrischen Funktionen
  und von Interpolationspolynomen.
\newblock {\em Numer. Math.\/}, {\bf vol.~2}:pp.~238--251, 1973.

\bibitem{BeWe93}
V.~Weispfenning and T.~Becker.
\newblock {\em Groebner bases: a computational approach to commutative
  algebra\/}, vol. 141 of {\em Graduate Texts in Mathematics: readings in
  mathematics\/}.
\newblock Springer, 1993.

\end{thebibliography}

}
\center{\small
\begin{tabular}{ll}
{\sc J. Heintz} 			& {\sc Guillermo Matera}\\
Depto. de Matem\'aticas, Est. y Comp.	& Depto. de Matem\'aticas\\
Facultad de Ciencias			& FCEyN, Universidad de Buenos Aires\\
Universidad de Cantabria 		& Ciudad Universitaria, Pab. I\\
E-39071 SANTANDER, Spain		& (1428) BUENOS AIRES, Argentina\\
{\tt heintz@matsun1.matesco.unican.es} 	& {\tt gmatera@dm.uba.ar}\\
Depto. de Matem\'aticas			& Instituto de Ciencias, Roca 850\\
FCEyN, Universidad de Buenos Aires	& Universidad Nacional de Gral. Sarmiento\\
Ciudad Universitaria, Pab. I		& (1663) San Miguel - Pcia. de Buenos Aires\\
(1428) BUENOS AIRES, Argentina		& Argentina\\
\\
{\sc Luis M. Pardo}			& {\sc Rosita Wachenchauzer}\\
Depto. de Matem\'aticas, Est. y Comp.	& Departamento de Computaci\'on\\
Facultad de Ciencias			& FCEyN, Universidad de Buenos Aires\\
Universidad de Cantabria		& Ciudad Universitaria\\
E-39071 SANTANDER, Spain		& (1428) BUENOS AIRES, Argentina\\
{\tt pardo@matsun1.matesco.unican.es}	& {\tt rosita@dc.uba.ar}\\
\end{tabular}
}


\end{document}
